\documentclass[11pt,draft]{article}
\usepackage{amsfonts}
\usepackage{amsfonts}
\usepackage{mathrsfs}
\usepackage{amsmath,amsfonts,mathrsfs,amssymb,color}
\usepackage{indentfirst}

\numberwithin{equation}{section}

\newtheorem{theo.}{\quad\, Theorem}[section]
\newtheorem{defi.}{\quad\, Definition}[section]
\newtheorem{lemm.}{\quad\, Lemma}[section]
\newtheorem{coro.}{\quad\, Corollary}[section]

\begin{document}

\title {Global structure of radial sign-changing solutions for the prescribed mean curvature problem in a ball }
\author{
Ruyun Ma$^{a,*}$\ \ \ \ \ \ \ Hongliang Gao$^{b}$
\\
{\small Department of Mathematics, Northwest
Normal University, Lanzhou 730070, P R China}\\
}
\date{} \maketitle
\noindent\footnote[0]{E-mail addresses: mary@nwnu.edu.cn (R.Ma), gaohongliang101@163.com (H.Gao)
} \footnote[0] {$^*$Supported by the
NSFC (No.11361054), SRFDP(No.20126203110004) and Gansu provincial National Science Foundation of China (No.1208RJZA258). }

{\small\bf Abstract.} {\small In this paper, we are concerned with the global structure of radial solutions, with prescribed nodal
properties, to the boundary value problem
 $$\text{div}\big(\phi_{N}(\nabla v)\big)+\lambda f(|x|, v)=0 ~~~\text{in} ~~B(R), ~~~ v=0 ~~~\text{on} ~~\partial B(R),
 $$
 where $\phi_{N}(y)=\frac{y}{\sqrt{1-|y|^{2}}},\; y\in \mathbb{R}^{N}$, $\lambda$ is a positive parameter, $B(R)=\{x\in \mathbb{R}^{N} :|x|<R\}$,  and $|\cdot|$ denote the Euclidean norm in $\mathbb{R}^{N}$.  All results, depending on the behavior of nonlinear term $f$ near 0, are obtained by using global bifurcation techniques.

\vskip 3mm


\vskip 3mm

{\small\bf MR(2000)\ \ \ 34B10, \ 34B18}

\baselineskip 19 pt

\section{Introduction}

    In this paper, we are concerned with the global structure of radial solutions, with prescribed nodal properties, to the boundary value problem
$$\text{div}\big(\phi_{N}(\nabla v)\big)+\lambda f(|x|, v)=0 ~~~\text{in} ~~B(R), ~~~ v=0 ~~~\text{on} ~~\partial B(R),
\eqno (1.1)
$$
 where $\phi_{N}(y)=\frac{y}{\sqrt{1-|y|^{2}}},\; y\in \mathbb{R}^{N}$, $\lambda$ is a positive parameter, $B(R)=\{x\in \mathbb{R}^{N} :|x|<R\}$,  and $|\cdot|$ denote the Euclidean norm in $\mathbb{R}^{N}$, $f$ satisfies

 \vskip 2mm

 \noindent$(H_1)$ $f: [0,R]\times (-\alpha,\alpha)\to \mathbb{R}$ is a continuous function, with $R<\alpha\leq \infty$ and such that  $f(r,s)s>0$ for $r\in [0,R]$  and $s\in (-\alpha,0)\cup(0,\alpha)$.

 \vskip 3mm

Dirichlet problem (1.1) is associated to mean curvature operator in flat Minkowski space
$$\mathbb{L}^{N+1}:=\{(x,t): x\in\mathbb{R}^{N}, t\in\mathbb{R}\}$$ endowed with the Lorentzian metric
$$\Sigma_{j=1}^{N}(dx_{j})^{2}-(dt)^{2},$$
where $(x,t)$ are the canonical coordinates in $\mathbb{R}^{N+1}$.

It is known (see e.g. [1, 4, 12, 29, 35]) that the study of spacelike submanifolds of codimension one in $\mathbb{L}^{N+1}$
with prescribed mean extrinsic curvature leads to Dirichlet problems of the type
$$\mathcal{M}v=H(x, v)\ \  \ \text{in}\ \ \ \ \ \Omega, \ \ \ \ \ \ v=0 \ \ \ \text{on}\ \partial \Omega,
\eqno (1.2)
$$
where
$$\mathcal{M}v=\text{div} \Big(\frac{\nabla v}{\sqrt{1-|\nabla v|^2}}\Big),
$$
$\Omega$ is a bounded domain in $\mathbb{R}^N$ and the nonlinearity $H: \Omega\times \mathbb{R}\to \mathbb{R}$ is continuous.

The starting point of this type of problems is the seminal paper [12] which deals with entire
solutions of $\mathcal{M}v= 0$. The equation $\mathcal{M}v=\text{constant}$ is then analyzed in
[35], while $\mathcal{M}v= f(v)$ with a general nonlinearity $f$ is considered in [9]. More general sign changing nonlinearities are studied in [5].

\vskip 3mm

If $H$ is bounded, then it has been shown by Bartnik and Simon [4] that (1.2) has at least one
solution $u \in C^1(\Omega)\cap W^{2,2}(\Omega)$.
Also, when $\Omega$ is a ball or an annulus in $\mathbb{R}^N$ and the nonlinearity
$H$ has a radial structure, then it has been proved in [6] that (1.2)
has at least one classical radial solution. This can be seen as a {\it universal} existence result for
the above problem in the radial case.

 Very recently, Bereanu, Jebelean and Torres [7] use Leray-Schauder degree arguments and critical point theory for convex, lower
semicontinuous perturbations of $C^1$-functionals, proved existence of classical positive radial solutions
for Dirichlet problems
$$\mathcal{M}v+f(|x|, v)=0 \ \ ~~~\text{in} ~~B(R), ~~~ v=0 ~~~\text{on} ~~\partial B(R),
$$
under the condition $(H_1)$ and
$$\lim_{s\to 0}\frac{f(r,s)}{\phi_1(s)}=\infty\ \ \  \ \text{uniformly for}\ r\in [0, R].
\eqno (1.3)
$$
Bereanu, Jebelean and Torres [8] used the upper and lower solutions and Leray-Schauder degree type arguments to study the special case of
$$\mathcal{M}v+\lambda \mu(|x|)v^q=0 \ \ ~~~\text{in} ~~B(R), ~~~ v=0 ~~~\text{on} ~~\partial B(R),
\eqno(1.4)
$$
under the condition

\noindent{($H_2$)} $N\geq 2$ is an integer, $R > 0,\; q > 1$ and $\mu: [0,\infty)\to \mathbb{R}$ is continuous, $\mu(r) > 0$ for all
$r >0$.

 \noindent They proved that there exists $\Lambda>0$ such that problem (1.4) has
zero, at least one or at least two positive radial solutions according to $\lambda\in (0,\Lambda)$, $\lambda = \Lambda$ or $\lambda > \Lambda$.
Moreover, $\Lambda$ is strictly decreasing with respect to $R$.

 \vskip 3mm

However,  the existence of sign-changing radial solutions for (1.1) has been scarcely explored in the related literature, see Capietto, Dambrosio and Zanolin [10].
\vskip 3mm

 When dealing with radial solutions to (1.1) on a ball, one is led to study (setting $|x| = r$)
the BVP
$$\left\{\aligned
 &(r^{N-1}\phi_1(u'))'+\lambda r^{N-1} f(r, u)=0, \ \ \ r\in (0, R),\\
 &u'(0)=0=u(R).\\
 \endaligned
 \right.
 \eqno (1.5)
 $$
 \vskip 3mm

Capietto, Dambrosio and Zanolin [10] used  a degree approach combined with a time-map technique (i.e. Abstract continuation theorem in [19]) to establish the following
\vskip 3mm

 \noindent{\bf Theorem A.} [10, Theorem 3.2] Assume that

\noindent$(H_f)$ $f(r, 0)\equiv 0$ and
$$\lim_{s\to 0}\frac{f(r,s)}{\phi_1(s)}=\infty, \ \ \ \ \text{uniformly in} \  r\in [0, R];$$

\noindent$(H_F)$  $F(r,s)=\int^s_0 f(r,x)dx$, and $F(r,s)$ is differentiable with respect to $r\in [0, R]$ and there exists a continuous
function $\omega: [0, R]\to \mathbb{R}^+$ such that
$$\Big|\frac{\partial F}{\partial r}(r,s)\Big|\leq \omega(r)F(r,s), \ \ \ \ r\in [0, R], \ s\in (-\alpha, \alpha).$$
  Then,  there exists $n_0$ such that for every $n >n_0$, $(1.5)_{\lambda=1}$ has at least two  solutions
$u^+_n$ and $v^-_n$ with $u^+_n(0) > 0$, $v^-_n(0) < 0$, all having exactly $n$
 zeros in $[0, R)$. Moreover, we have
$$\lim_{n\to\infty}\; |u^+_n(r)|+|(u_n^+)'(r)|=0= \lim_{n\to\infty}\; |v^-_n(r)|+|(v_n^-)'(r)| \ \ \ \ \text{uniformly in}\ r\in [0,R].
$$

 \vskip 3mm

 Motivated above papers, in this paper, we investigate the global structure of radial solutions, with prescribed nodal properties, to Dirichlet problem (1.1)
by the unilateral global bifurcation theory of Dancer [16-17] and L\'{o}pez-G\'{o}mez [24, Sections 6.4, 6.5] and some preliminary results on the superior limit of a sequence of connected components due to  Ma and An [26-27]. We shall use the following assumptions

\vskip 2mm

\noindent(A1) $R\in (0, \infty)$ and $\delta\in [0, R)$,  $f: [0,R]\times (-\alpha,\alpha)\to \mathbb{R}$ is a continuous function, with $R<\alpha\leq \infty$ and such that  $f(r,s)s>0$ for $r\in [0,R]$  and $s\in (-\alpha,0)\cup(0,\alpha)$.

\vskip 2mm

\noindent(A2) $\lim\limits_{s\to 0}\frac{f(r,s)}{s}=m(r)$ uniformly $r\in [\delta,R]$ for some $m\in C[\delta,R]$ and
$$m(r)\geq 0, \ \ m(r)\not\equiv 0\ \ \text{on any subinterval of}\ [\delta,R];$$

\vskip 2mm

\noindent(A3) $f(r, 0)\equiv 0$ and
$$f_0:=\lim_{s\to 0}\frac{f(r,s)}{\phi_1(s)}=\infty, \ \ \ \ \text{uniformly in} \  r\in [\delta, R];$$

\vskip 2mm

\noindent(A4)   $F(r,s)$ is differentiable with respect to $r\in [\delta, R]$ and there exists a continuous
function $\omega: [\delta, R]\to \mathbb{R}^+$ such that
$$\Big|\frac{\partial F}{\partial r}(r,s)\Big|\leq \omega(r)F(r,s), \ \ \ \ r\in [\delta, R], \ s\in (-\alpha, \alpha).
$$

\vskip3mm

 To study the global structure of radial solutions of problem (1.1), we need to study the family of auxiliary problems
 $$\left\{\aligned
 &(r^{N-1}\phi_1(u'))'+\lambda r^{N-1} f(r, u)=0,\ \ \ \ r\in (\delta, R),\\
 &u'(\delta)=0=u(R).\\
 \endaligned
 \right.
 \eqno (1.6)_\delta
 $$
\vskip 3mm

For given $\delta\in [0, R)$. Let
$$X_\delta=C[\delta, R],\ \ \ \ \ \ E_\delta=\{u\in C^1[\delta, R]:\, u'(\delta)=u(R)=0\}$$
be the Banach spaces endowed with the normals
$$||u||_{C[\delta, R]}=\sup_{r\in [\delta, R]}|u(r)|, \ \ \ \ \    ||u||_{C^1[\delta, R]}=\sup_{r\in [\delta, R]}|u(r)|+\sup_{r\in [\delta, R]}|u'(r)|,
$$
respectively. Denoted by  $\Sigma_\delta$ be the closure of the set
$$\{(\lambda, u)\in [0, \infty)\times C^1[\delta, R]:
 u\ \text{satisfies} \ (1.6)_\delta, \ \text{and}\ u\not\equiv 0\}
$$
in $\mathbb{R}\times E_\delta$.
Let $S^+_{k,\delta}$ denote the set of function $u\in E_\delta$,
which have exactly $k-1$ non-degenerate nodal zeros in $(\delta, R)$ and there exists $\sigma_1>0$ such that $u$ is positive in $(\delta, \delta+\sigma_1)$, and set $S_{k,\delta}^-=-S_{k,\delta}^+$, and $S_{k,\delta} = S^+_{k,\delta}\cup  S^-_{k,\delta}$. It is clear that $S^+_{k,\delta}$ and $S^-_{k,\delta}$ are disjoint and open in $E_\delta$.
 Finally, let $\Phi^\pm_{k,\delta} = \mathbb{R}\times S^\pm_{k,\delta}$ and $\Phi_{k,\delta} = \mathbb{R}\times S_{k,\delta}$ under the product topology.
Denoted by $\theta$ be the zero element in $E_\delta$.
\vskip 3mm

Let $\delta\in [0, R)$ be given and let $\lambda_k(m, \delta)$ be the $k$-th eigenvalue of
$$\left\{\aligned
&-(r^{N-1}u')'=\lambda r^{N-1}m(r)u,\ \ \ \ r\in (\delta, R),\\
&u'(\delta)=u(R)=0.\\
\endaligned
\right.
\eqno (1.7)_\delta
$$

The main results of the paper are as follows

\vskip 3mm

\noindent{\bf Theorem 1.1.}\
 Assume that (A1) and (A2) hold. Then for $\nu\in\{+, -\}$ and  $k\in \mathbb{N}$, there exists a connected component $\zeta^\nu\in \Sigma_\delta$, such that

(a) $\big(\zeta^\nu\setminus\{(\lambda_k(m, \delta), \theta)\}\big)\subset \big((0,\infty)\times  \text{int}\, \Phi^\nu_{k,\delta}$\big);

(b) $\zeta^\nu$ joins $(\lambda_k(m, \delta),\theta)$ with infinity in $\lambda$ direction;

(c) $\text{Proj}_\mathbb{R}\,\zeta^\nu=[\lambda_*, \infty)\subset (0, \infty)$ for some $\lambda_*>0$.

\vskip 3mm

\noindent{\bf Theorem 1.2.}\ Let $\delta\in [0, R)$ be given. Assume that (A1), (A3) and (A4) hold. Then for $\nu\in\{+, -\}$ and  $k\in \mathbb{N}$, there exists a connected component $\zeta^\nu\in \Sigma_\delta$ such that

(a) $\big(\zeta^\nu\setminus\{(0, \theta)\}\big)\subset\big((0,\infty)\times  \text{int}\, \Phi^\nu_{k,\delta}\big)$;

(b) $\zeta^\nu$ joins $(0,\theta)$ with infinity in $\lambda$ direction;

(c) $\text{Proj}_\mathbb{R}\,\zeta^\nu=[0, \infty)$.

\vskip 3mm

\vskip 3mm

Obviously, as the immediate consequences of Theorem 1.1-1.2, we have the following

\vskip 3mm

\noindent{\bf Corollary 1.1.}\ Let $\delta\in [0, R)$ be given. Assume that (A1) and (A2) hold. Then for $\nu\in\{+, -\}$ and  $k\in \mathbb{N}$, there exists $\lambda_{\ast}^\nu\in (0,\lambda_k(m,\delta)]$ such that, for all $\lambda\in(0,\lambda_{\ast}^\nu)$, the problem
$(1.6)_\delta$ has no solution in $S^\nu_{k,\delta}$ and, for all $\lambda>\lambda_k(m,\delta)$ has at least one solution in $S^\nu_{k,\delta}$.

\vskip 3mm

\noindent{\bf Corollary 1.2.}\ Let $\delta\in [0, R)$ be given. Assume that (A1), (A3) and (A4) hold. Then for $\nu\in\{+, -\}$ and  $n\in \mathbb{N}$, the problem
$(1.6)_\delta$ has at least one solution $u_n^\nu\in S^\nu_{n,\delta}$ for $\lambda=1$. Moreover,
$$\lim_{n\to\infty}\; |u^\nu_n(r)|+|(u_n^\nu)'(r)|=0 \ \ \text{uniformly in}\ r\in [0,R].
$$

\vskip 3mm

\noindent{\bf Remark 1.1.}\ Corollary 1.2 guarantees that for each  $n\in \mathbb{N}$, $(1.6)_\delta$ has two solutions
$u^+_n$ and $v^-_n$ with $u^+_n(0) > 0$, $v^-_n(0) < 0$, for all $\lambda>0$ and $\delta\in [0, R)$.
[11, Theorem 3.2] deal with the more general problem
$$\left\{\aligned
 &\text{div} \big(a(|\nabla u|)\nabla u\big)+f(|x|, u)=0 ~~~\text{in} ~~B(R),\\
 &u=0 ~~~\text{on} ~~\partial B(R),\\
 \endaligned
 \right.
  \eqno (1.8)
  $$
 where $a:[0,\epsilon_1]\to [0, +\infty),\; (\epsilon_1>0)$. However, they only proved that there exists $n_0\in \mathbb{N}$ such that for $n\geq n_0$, (1.8) has at least two radial solutions $u^+_n$ and $v^-_n$  with $u^+_n(0) > 0$, $v^-_n(0) < 0$, all having exactly $n$
 zeros in $[0, R)$.

\vskip 3mm

\noindent{\bf Remark 1.2.} Notice that $(1.6)_0$ a singularity appears for $r = 0$. Beside this intrinsic aspect of  $(1.6)_0$, assumption (A3) means that $\frac{\partial f}{\partial s}$ does not exist at $s=0$.
Hence, in particular, when one tries to develop some shooting argument, global
existence and uniqueness to initial value problems associated to the equation in $(1.6)_\delta$ are not
guaranteed. This is one of the reasons why few results for the case $f_0=\infty$
 is available in the literature:  we refer to the earlier works of Krasnosel'skii, Perov, Povolotskii and Zabreiko [22, Section 15], Shekhter [34, Section 15] for a more classical approach in the ODE's case. For other results, in the PDE's setting, we also refer to  Omari and Zanolin [31],  Willem [37], mainly for the
case $a(\cdot)\equiv 1$.

\vskip 3mm

        The rest of the paper is organized as follows. In Section 2 we state some preliminary results on the superior limit of a sequence of connected components due to Ma and An [26-27] and on the Unilateral global bifurcation theorem of Dancer [16-17] and L\'{o}pez-G\'{o}mez [24]. Section 3 is devoted to establish the existence of connected component of radial solutions for the prescribed mean curvature problem in an annular domain via global bifurcation technique. Finally in Section 4, we shall use the components obtained in Section 3 to construct the desired components of radial solutions for the prescribed mean curvature problem in a ball and prove Theorem 1.1-1.2 and their corollaries.

\vskip 3mm

        For other results concerning the problem associated to prescribed mean curvature
equations in Minkowski space we refer the reader to [5, 9, 23, 29].
   The existence of radial solutions satisfying various boundary conditions has been investigated
by many authors, see
 Esteban [18], Castro and Kurepa
[11], Grillakis [20], Guo [21],
 Cheng [13], Ambrosetti, Garcia Azorero and Peral [2], Njoku, Omari and Zanolin
[30],  Dai and Ma [15] and references therein.

\vskip 3mm

\section{Preliminary results}

 \subsection{Unbounded connected component}

Let $M$ be a metric space and $\{C_n\,|\, n=1,2,\cdots\}$  a family
of subsets of $M$. Then the superior limit $\mathscr{D}$ of
$\{C_n\}$ is defined by
$$\mathscr{D}:=\limsup\limits_{n\to\infty}C_n=\{x\in M\,|\,\exists\
\{n_k\}\subset \mathbb{N},\ x_{n_k}\in C_{n_k},\ \text{such that}\
x_{n_k}\to x\}. \eqno(2.1)$$

Let $X$ be a Banach space with the norm $\|\cdot\|$. A {\it component} of a set $M\subset X$ means a maximal connected subset of $M$, see [36] for the detail.

\vskip 3mm

The following results are somewhat scattered in Ma and An [26-27], Ma and Gao [28].
\vskip 2mm

\noindent{\bf Lemma 2.1 ([26, Lemma 2.4; 27, Lemma 2.2])}\ Let $X$
be a Banach space. Let $\{C_n\}$ be a family of closed connected
subsets of $X$. Assume that

(i) there exist $z_n\in C_n,\, n=1,2,\cdots$ and $z_\ast\in X$ such
that $z_n\to z_\ast$;

(ii)
$\lim\limits_{n\to\infty}r_n=\lim\limits_{n\to\infty}\sup\{\|u\|\,|\,u\in
C_n\}=\infty$;

(iii) for every $R>0$, $(\bigcup_{n=1}^\infty C_n)\cap B_R$ is a
relatively compact of $X$. \\
Then there exists an unbounded component $\mathscr{C}$ in
$\mathscr{D}$ and $z_\ast\in \mathscr{C}$.

\vskip 3mm

\subsection{Global alternative of Rabinowitz}

In order to formulate and prove main results of this section, it is convenient to introduce Dancer [16-17]
and L\'{o}pez-G\'{o}mez's notations [24].
Let
$\mathbb{X} = \mathbb{R}\times X$.
Given any $\mu\in \mathbb{R}$ and $0 < s <+\infty$, we consider an open neighborhood
of $(\mu, 0)$ in $\mathbb{X}$ defined by
$$\mathbb{B}(\mu, 0):=\{(\mu,u)\in \mathbb{X}: ||u||+|\mu|<s\}.
$$

A mapping $G:\mathbb{X}\to X$ is said to satisfy {\it Assumption $\mathfrak{A}$} if

\noindent(1) $G(0, \lambda)=0$ for $\lambda\in \mathbb{R}$;

\noindent(2) $G$ is completely continuous, $G(x, \lambda)=\lambda Lx+H(\lambda, x)$,
where $L$ is a continuous linear operator on $X$;

\noindent(3) $||H(\lambda, x)||/||x||\to 0$ as $||x||\to 0$ uniformly on bounded subsets of $\mathbb{R}$.

\vskip 3mm

Define $\Phi:\mathbb{X}\to X$ by
$$\Phi(\lambda, x)=x-G(\lambda, x)$$
and
 $$\mathcal{S}:=\overline{\{(\mu, u)\in \mathbb{X}:\ \Phi(\mu, u)=0, \ u\neq 0\}}^\mathbb{X}.
 $$

Assume that $\mu\in r(L)$ such that $\mu$ has algebraic multiplicity $1$.
Suppose that $\varphi\in X\setminus\{0\}$ such that
$$\varphi=\mu L\varphi.$$
Let $X_0$ be a closed subspace of $X$ such that
$$X=\text{span}\{\varphi\}\oplus X_0.
$$

Let
$C_\mu$ to be the component of $\mathcal{S}$ containing $(\mu, 0)$.

\vskip 3mm

\noindent{\bf Theorem B.  Global bifurcation of Rabinowitz, see L\'{o}pez-G\'{o}mez [24, Corollary 6.3.1]} \\
\ Assume that $\mu\in r(L)$ has algebraic multiplicity $1$.  Then, one of the following non-excluding options occurs. Either
\\
1. $C_\mu$ is unbounded in $\mathbb{R}\times X$.
\\
2. There exists $\mu_1\in r(L)\setminus\{\mu\}$ such that $(\mu_1, 0)\in C_\mu$.\hfill{$\Box$}

\subsection{Unilateral global bifurcation theorem}

\vskip 3mm

In this subsection, we shall introduce the unilateral global bifurcation theorem,
see Dancer [16-17] and
  L\'{o}pez-G\'{o}mez [24].

\vskip 2mm

According to the Hahn-Banach theorem, there exists a linear functional $l\in X^*$, here $X^*$ denotes the
dual space of $X$, such that
$$l(\varphi)=1, \ \ \ X_0=\{u\in X:\, l(u)=0\}.
$$
Finally, for any $0 < \eta < 1$, we define
$$K_\eta:=\{(\mu, u)\in \mathbb{X}:\, |l(u)|>\eta ||u||\}.
$$
Since
$$u\mapsto |l(u)|-\eta ||u||$$
is continuous, $K_\eta$ is an open subset of $\mathbb{X}$ consisting of two disjoint components $K^+_\eta$ and $K^-_\eta$, where
$$K^+_\eta:=\{(\mu, u)\in \mathbb{X}:\, l(u)>\eta ||u||\}.\ \
$$
$$K^-_\eta:=\{(\mu, u)\in \mathbb{X}:\, l(u)<-\eta ||u||\}.
$$
In particular, both $K^+_\eta$ and $K^-_\eta$ are convex cones, $K^+_\eta=-K^-_\eta$, and $\nu t\varphi \in K^\nu_\eta$ for every $t > 0$, where
$\nu\in \{+, -\}$.
Applying the similar method to prove [24, Lemma 6.4.1] with obvious changes, we may obtain the
following result.

\vskip 3mm

\noindent{\bf Lemma 2.2.} For every $\eta\in (0, 1)$, there exists a number $\delta_0> 0$ such that for each $0 <\delta<\delta_0$,
$$((S\setminus\{(\mu, 0)\})\cap\bar B_\delta(\mu, 0))\subset K_\eta.
$$
Moreover, for each
$$(\mu, u)\in (S\setminus\{(\mu, 0)\})\cap\bar B_\delta(\mu, 0)),
$$
there are $s\in \mathbb{R}$ and a unique $y \in X_0$ such that
$$u=s\varphi+y\ \ \ \text{and}\ \ |s|>\eta ||u||.
$$
Furthermore, for these solutions $(\lambda, u)$
$$\lambda=\mu+\circ (1)\ \ \ \text{and}\ \  y=\circ (s)$$
as $s\to 0$.

\vskip 3mm

Let $\delta > 0$ be the constant from Lemma 2.2. For $0 <\epsilon<\delta$ we define
$\mathcal{D}^\nu_{\mu, \epsilon}$ to be the component
of $\{(\mu, \theta)\} \cup (\mathcal{S} \cap \bar B_\epsilon \cap K^\nu_\eta)$ containing $(\mu, \theta)$, $\mathcal{C}^\nu_{\mu, \epsilon}$ to be the component of $C_\mu\setminus \mathcal{D}^{-\nu}_{\mu, \epsilon}$ containing
$(\mu, \theta)$, and $C^\nu_\mu$ to be the closure of $\cup_{0<\epsilon\leq \delta} C^\nu_{\mu, \epsilon}$. Clearly, $C^\nu_\mu$ is connected. Thanks to Lemma 2.2, the
definition of $C^\nu_\mu$ is independent from the choice of $\eta$ and
$$C_\mu = C^+_\mu\cup C^-_\mu.
$$

\noindent{\bf Theorem C. \  Unilateral global bifurcation, see Dancer [16-17]}
\\
Either $C^+_\mu$ and $C^-_\mu$ are both unbounded or $C^+_\mu \cap C^-_\mu\neq \{(\mu, \theta)\}$.\hfill{$\Box$}

\vskip 3mm

\subsection{Uniqueness of solutions of Cauchy problem}

\vskip 3mm

\noindent{\bf Lemma 2.3.} Let $\tau\in (\delta, R)$ be given. Let
$f:[\delta, R]\times \mathbb{R}\to \mathbb{R}$ be continuous and  let $f(r, u)$ be Lipschitz continuous in $u$ on bounded sets, and $f(r,0)\equiv 0$ for $r\in [\delta, R]$.
If $u$ is a solution of
$$
\left\{\aligned
&(r^{N-1}\phi(u'))'+\lambda r^{N-1}f(r, u)=0,\\
&u'(\tau) = 0=u(\tau),\ \ \  \\
\endaligned
\right.
\eqno (2.2)$$
then $u\equiv 0$ in $[\delta, R]$.
\vskip 3mm

\noindent{\bf Proof.} (2.2) is equivalent to
$$
\left\{\aligned
&-(r^{N-1}u')'=r^{N-1}[\lambda f(r,u)h(u')-\frac{N-1}ru'^{3}],\ \ \  \ r\in (\delta, R),\\
&u'(\tau) = 0=u(\tau).  \\
\endaligned
\right.
$$
Since the function
$$g(r, u, p):=r^{N-1}[\lambda f(r,u)(1-p^2)^{3/2}-\frac{N-1}rp^{3}]$$
is continuous in $[\delta, R]\times \mathbb{R}\times [-1, 1]$ and is Lipschitz in $(u,p)$ on any bounded subset of $\mathbb{R}\times [-1, 1]$, it deduce that $u\equiv 0$ in $[\delta, R]$.  \hfill{$\Box$}

\vskip 3mm

Finally, we consider the case that $f$ does not meet the Lipschitz condition at $u=0$, i.e. we allow $f_0=\infty$.

\vskip 3mm

\noindent{\bf Lemma 2.4.} [10, Lemma 2.3] Let $(A4)$ hold.
Let $\epsilon_0=0.9$. Then, for every $\epsilon\leq \epsilon_0$ there exists $d_\epsilon\in (0, \epsilon]$ such that if $u$ is a (local) solution of
$$\left\{\aligned
&(r^{N-1}\phi_1(u'))'+\lambda r^{N-1}f(r, u)=0,\\
&u(\delta)= d,\ \ \  u'(\delta) = 0\\
\endaligned
\right.
$$
with $|d|\leq d_\epsilon$, then $u$ can be defined on $[\delta, R]$ and $||u||_{C^1[\delta, R]}\leq \epsilon$.

\vskip 3mm

\noindent{\bf Lemma 2.5.} [10, Lemma 2.5] Let $(A4)$ hold. Let $u$ be a solution of
$$\left\{\aligned
&(r^{N-1}\phi_1(u'))'+\lambda r^{N-1}f(r, u)=0,\\
&u(r_0)= u'(r_0) = 0, \ \ \ \ r_0\in (\delta, R].\\
\endaligned
\right.
$$
Then $u\equiv 0$ in $[\delta, R]$.
\hfill{$\Box$}

\vskip 3mm

\section{Radial solutions for the prescribed mean curvature problem in an annular domain}

Let $\delta\in (0,R)$ be a given constant in this section.

\vskip 3mm

Let us consider the following boundary value problem
$$
\left\{\aligned
&\text{div}\big(\phi_{N}(\nabla v)\big)+\lambda f(|x|, v)=0 ~~~\text{in} ~~\mathcal{A},\\
&\frac{\partial v}{\partial \nu}=0 ~~~\text{on} ~~\Gamma_{1},  \ \ \ \ \ \ \ \ \ \  v=0 ~~~\text{on} ~~\Gamma_{2}, \\
\endaligned
\right.
\eqno(3.1)
$$
where $$\mathcal{A}=\{x\in \mathbb{R}^{N} :\delta<|x|<R\},$$ $$\Gamma_{1}=\{x\in \mathbb{R}^{N} :|x|=\delta\}, \ \ \ \ \ ~\Gamma_{2}=\{x\in \mathbb{R}^{N} :|x|=R\},$$
 $\frac{\partial v}{\partial \nu}$ and $|\cdot|$ denote the outward normal derivative of $v$ and the Euclidean norm in $\mathbb{R}^{N}$, respectively.

 \vskip 3mm

Setting, as usual, $|x|=r$ and $v(x)=u(r)$, the above problem (3.1) reduces to
$$
\left\{\aligned
&-(r^{N-1}\phi_{1}(u'))'=\lambda r^{N-1}f(r,u),\\
&u'(\delta)=0=u(R).\\
\endaligned
\right.
\eqno (3.2)_\delta
$$
To find a radial solution of (3.1),
it is enough to find a solution of $(3.2)_\delta$.

\vskip 3mm

\noindent{\bf Remark 3.1.}\ It is worth remarking that $(3.2)_\delta$
is equivalent to
$$\left\{\begin{array}{ll}
-(r^{N-1}u')'= r^{N-1}[\lambda f(r,u)h(u')-\frac{N-1}ru'^{3}],\ \ \  \ r\in (\delta, R),\\
u'(\delta)=0=u(R).  \\
\end{array}
\right. \eqno(3.3)_\delta
$$
Since the nonlinearity $F(r,u,p):=\lambda f(r,u)h(p)-\frac{N-1}rp^{3}$ is
singular at $r=0$ if $\delta=0$, we cannot deal with $(3.3)_0$ via
the spectrum of $(1.4)_0$ directly. However, $F(r,u,p)$ is regular
at $r=\delta$ if $\delta>0$, in this case, $(3.3)_\delta$ with
$\delta>0$ can be treated via the spectrum of $(1.4)_\delta$ and the
standard bifurcation technique. This is why we firstly study the
prescribed mean curvature problem in an annular domain.

 \vskip 3mm

\noindent{\bf Lemma 3.1.}  Let $u\in S^+_{k,\delta}$ be a solution of
$$
(r^{N-1}\phi_1(u'))'+r^{N-1}f(r, u)=0,\ \ \ r\in (\delta, R)
\eqno (3.4)
$$
Assume that (A1) holds.
Then
$u(\delta)>0$.

\vskip 2mm

\noindent{\bf Proof.} Denote by $\tau_1$ the first positive zero of $u$.
Let $u(\delta):=d$.
Then $u$ satisfies
$$
\left\{\aligned
&(r^{N-1}\phi_1(u'))'+r^{N-1}f(r, u)=0,\ \ \ r\in (\delta,\tau_1),\\
&u'(\delta)=0, \ \ \ \ \ u(\tau_1)=0.\\
\endaligned
\right.
\eqno (3.5)
$$
It is easy to check that (3.5) is equivalent to
$$r^{N-1}\phi_1(u'(r))=-\int^r_\delta t^{N-1}f(t, u(t))dt, \ \ \ \ r\in [\delta, \tau_1]
$$
it follows $u'\leq 0$ because $f(r,s)\geq 0$ for all $r\in [\delta,R]$ and $s\in[0,\alpha)$, so $u$ is decreasing. Since
$u(\tau_1) = 0$, we have $u\geq 0$ on $[\delta,\tau_1]$. As $u$ is not identically zero, one has $u(\delta) > 0$,  and subsequently, $u'
 < 0$ on $(\delta,\tau_1]$, which ensures that actually $u$ is strictly decreasing in $[\delta, \tau_1]$ and $u > 0$
on $[\delta,\tau_1)$.
\hfill{$\Box$}

\vskip 3mm

 \vskip 3mm

 \noindent{\bf Definition 3.1.} For $y\in S_k$, denote the zeros of $y$ by
 $$(\delta<)\,\tau_1<\cdots<\tau_k\,(=R).$$
If  $y'(\delta)=0$ and,  for each $j\in\{1, \cdots ,k-1\}$,
 there exists exactly one $\xi_j\in (\tau_{j}, \tau_{j+1})$, such that $y'(\xi_j)=0$,
  then we call that $y$ possesses the {\it Property $[P(k)]$}.

\vskip 3mm

 Let
 $$\Pi_k=\{y\in C^1[\delta, R]:\, y'(\delta)=y(R)=0,\, y\ \text{possesses the Property $[P(k)]$} \}.$$

\vskip 3mm

\vskip 3mm

\subsection{Eigenvalue problem in an annular domain}

Let $\delta\in (0, R)$ be given.  Let us recall the weighted eigenvalue problem
$$\left\{\begin{array}{ll}
-(r^{N-1}u')'=\lambda r^{N-1}m(r)u,\ \ \  \ r\in (\delta, R),\\
u'(\delta)=0=u(R),  \\
\end{array}
\right. \eqno(3.7)_\delta
$$
where

\noindent{(A5)} $m\in C[\delta, R]$ and $m(r)\geq 0, m(r)\not\equiv 0$ on any subinterval of $[\delta, R]$.

The following result is a special case of [32, Theorem 1.5.3] when $p=2$.

\noindent\textbf{Lemma~3.2.}\  \ Let (A5) hold. Then the eigenvalue problem $(3.7)_\delta$ has infinitely many simple real eigenvalues
$$0<\lambda_{1}(m, \delta)<\lambda_{2}(m, \delta)<\cdot\cdot\cdot<\lambda_{k}(m, \delta)<\cdots \to +\infty \ \ \text{as}\ \  k\to +\infty$$
and no other eigenvalues. Moreover, the algebraic multiplicity of $\lambda_k(m,\delta)$ is $1$, and the eigenfunction $\varphi_{k}$ corresponding to $\lambda_{k}(m, \delta)$ has exactly $k-1$ simple zeros in $(\delta, R)$.

\vskip 3mm

Let us  consider the following auxiliary problem
$$\left\{\begin{array}{ll}
-(r^{N-1}u')'=r^{N-1}h(r),\ \ \  \ \ r\in (\delta, R)\ \text{with}\ \delta>0,\\
u'(\delta)=0=u(R)  \\
\end{array}
\right.
$$
for a given $h\in X_\delta$. Its Green function for $N\geq 3$
is explicitly given by
$$K_\delta(t,s)=
\left\{
  \aligned
  &\frac 1{2-N}[R^{2-N}- t^{2-N}], \ \ \ \ \ \delta\leq s\leq t\leq R,\\
  &\frac 1{2-N}[R^{2-N}- s^{2-N}], \ \ \ \ \ \delta\leq t\leq s\leq R.\\
  \endaligned
  \right.
$$
and its Green function for $N=2$ is explicitly given by
$$K_\delta(t,s)
=
\left\{
  \aligned
  &\ln \frac Rt, \ \ \ \ \qquad\ \delta\leq s\leq t\leq R,\\
  &\ln \frac Rs, \ \ \ \ \qquad\ \delta\leq t\leq s\leq R.\\
  \endaligned
  \right.
$$
It is well-known that for every $h\in X_\delta$, the above auxiliary problem has a unique solution
$$u=\int^R_\delta K_\delta(t,s)s^{N-1}h(s)ds=:\mathcal{G}_\delta(h)
$$
It is easy to check that $\mathcal{G}_\delta : X_\delta \to E_\delta$ is continuous and compact (see [3]).

\vskip 3mm

Define a linear operator $\mathcal{L}_\delta: X_\delta\to E_\delta\ (\hookrightarrow X_\delta)$.
$$\mathcal{L}_\delta(u)(r):=\mathcal{G}_{\delta}(mu)(r).$$
Then $\mathcal{L}_\delta$ is compact and
$(3.7)_\delta$ is equivalent to
$$u=\lambda \mathcal{L}_\delta(u).
$$
Moreover,  $\mathcal{L}_\delta|_{E_\delta}: \ E_\delta\to E_\delta$ is compact.

\vskip 3mm

\subsection{An equivalent formulation}

 We will use some idea in Coelho et. al. [14].  Let us define a function $\tilde{f}: [\delta,R]\times \mathbb{R}\to \mathbb{R}$ by setting, for $r\in [\delta, R]$,
$$
\tilde{f}(r,s) =\left\{ \aligned
&f(r,s), \ \ \ \ \ \ \  \ \text{if}\ \ 0\leq |s|\leq R-\delta,\\
&0, \ \ \ \ \ \  \ \ \ \ \ \ \ \ \ \text{if}\ \ |s|\geq (R-\delta)+1,\\
&\text{linear}, \ \  \ \ \ \ \ \ \ \text{if}\ \ R-\delta<|s|<(R-\delta)+1.\\
\endaligned
\right.
$$
Since any solution $u$ of $(3.2)_\delta$ satisfies
 $$||u'||_{C[\delta,R]}<1, \ \ \ \ ||u||_{C[\delta,R]}< R-\delta,
 $$
 it follows that $(3.2)_\delta$ is equivalent to the same problem with $f$ replaced by $\tilde{f}$.  Clearly, $\tilde{f}$ satisfies all the properties assumed in the statement of the theorem. In the sequel, we shall replace $f$ with $\tilde{f}$; however, for the sake of simplicity, the modified function $\tilde{f}$ will still be denoted by $f$.
Next, let us define $h: \mathbb{R}\to \mathbb{R}$ by setting
$$
 h(y) =\left\{
 \aligned
&(1-y^{2})^{\frac{3}{2}}, \ \   \ \ \text{if}\ \ |y|\leq1,\\
&0,\ \ \ \ \ \  \ \ \ \ \ \ \ \ \ \text{if}\ \  |y|>1.\\
\endaligned
\right.
$$
{\it Claim.} A function $u\in C^{1}[\delta, R]$ is a solution of $(3.2)_\delta$ if and only if it is a  solution of the problem
$$\left\{\begin{array}{ll}
-(r^{N-1}u')'=\lambda r^{N-1}f(r,u)h(u')-(N-1)r^{N-2}u'^{3},\ \ \  \ r\in (\delta, R),\\
u'(\delta)=0=u(R).  \\
\end{array}
\right. \eqno(3.8)_\delta
$$
It is clear that a solution $u\in C^{1}[\delta, R]$ of  $(3.2)_\delta$  is a  solution of  $(3.8)_\delta$  as well. Conversely, suppose that $u\in C^{1}[\delta, R]$ is a  solution of  $(3.8)_\delta$ . We aim to show that
$$\|u'\|_{C[\delta,R]}<1.
$$
Assume by contradiction that this is not the case. Then we can easily find an interval $[a,b]\subseteq[\delta, R]$ such that, either $u'(a)=0,\ 0<|u'(r)|<1$ in $(a,b)$ and $|u'(b)|=1$, or $|u'(a)|=1,\ 0<|u'(r)|<1$ in $(a,b)$ and $u'(b)=0$. Suppose the former case occurs (in the latter one the argument would be similar). The function $u$ satisfies the equation
$$-\big(r^{N-1}\frac{u'}{\sqrt{1-u'^{2}}}\big)'=\lambda r^{N-1}f(r,u)$$
in $[a,b)$. For each $r\in (a,b)$, integrating over the interval $[a,r]$ and using (A1), we obtain
$$|\phi_{1}(u'(r))|=\Big|\frac{1}{r^{N-1}}\int_{a}^{r}\lambda t^{N-1}f(t,u)dt\Big|\leq M$$
and hence
$$|u'(r)|\leq\phi_{1}^{-1}(M)
$$
for every $r\in [a,b)$. Since $\phi_{1}^{-1}(M)<1$, taking the limit as $r\to b^{-}$ we get $|u'(b)|<1$. This is a contradiction.  Therefore $\|u'\|_{C[\delta,R]}<1$ and, as a consequence, $u$ is a  solution of $(3.2)_\delta$.

\subsection{Proof of Theorem 1.1-1.2 with $\delta\in (0, R)$}

In this subsection, we shall prove Theorem 1.1-1.2 in the case  $\delta>0$.
\vskip 3mm

\noindent{\bf Proof of Theorem 1.1 with $\delta\in (0, R)$.}\  \ By (A1) and (A2) we can write, for any $r\in [\delta,R]$ and every $s\in \mathbb{R}$,
$$f(r,s)=(m(r)+l(r,s))s,$$
where $l :[\delta,R]\times \mathbb{R}\to \mathbb{R}$ is a continuous function and
$$\lim_{s\to 0}l(r,s)=0
\eqno (3.9)
$$
uniformly in $[\delta,R]$. Let us set, for convenience, $k(y)=h(y)-1$ for $y\in \mathbb{R}$. We have
$$\lim_{y\to 0}\frac{k(y)}{y}=0.
\eqno(3.10)
$$
Define the operator $\mathcal{H}_\delta: \mathbb{R}\times E_\delta\to E_\delta$ by
$$\mathcal{H}_\delta(\lambda,u)(\cdot)=\mathcal{G}_\delta\big(\lambda[l(\cdot,u)+(m(\cdot)+l(\cdot,u))k(u')]u-\gamma(\cdot)u'^{3}\big)$$
where $\gamma(r)=\frac{N-1}{r}$.
Clearly, $\mathcal{H}_\delta$ is completely continuous and, by (3.9) and (3.10),
$$\lim_{\|u\|_{C^1[\delta,R]}\to 0}\frac{\|\mathcal{H}_\delta(\lambda,u)\|_{C^1[\delta,R]}}{\|u\|_{C^1[\delta,R]}}=0,
\eqno(3.11)
$$
uniformly with respect to $\lambda$ varying in bounded intervals. Observe that, for any $\lambda$, the couple $(\lambda,u)\in \mathbb{R}\times E_\delta$ is a solution of the equation
$$u=\lambda \mathcal{L}_\delta(u)+\mathcal{H}_\delta(\lambda,u)
\eqno(3.12)
$$
if and only if $u$ is a solution of $(3.2)_\delta$.

\vskip 2mm

Recall that  $\Sigma_\delta\subset\mathbb{R}\times E_\delta$ be the closure
of the set of all nontrivial solutions $(\lambda, u)$ of (3.12) with
$\lambda>0$. Note that the set $\{u\in E_\delta\,|\,
(\lambda,u)\in \Sigma_\delta\}$ is bounded in $E_\delta$.

\vskip 3mm

    As the algebraic multiplicity of $ \lambda_k(m, \delta)$  equals 1 [25], the local
index of $0$ as a fixed point of $\lambda \mathcal{L}_\delta$ changes sign as $\lambda$ crosses $ \lambda_k(m, \delta)$.
Therefore, according to a revised version of [24, Theorem 6.2.1],  there exists a component, denoted by $\mathfrak{C}_k\subset \Sigma_\delta$, emanating from $( \lambda_k(m, \delta), \theta)$.

\vskip 3mm

Now, we use some notations and preliminary results on Unilateral global bifurcation.

We shall show that  both  $\mathfrak{C}^+_k$ and $\mathfrak{C}_k^-$ are unbounded , and $\mathfrak{C}_k^\nu\subset \Phi^\nu_k$ for $\nu\in\{+,-\}$.

It is easy to check that
(3.12) enjoys the structural requirements for applying the unilateral global bifurcation theory of Dancer and L\'{o}pez-G\'{o}mez (see Theorem C} (by a counterexample of Dancer [16], the global unilateral theorem of Rabinowitz [33] is false as stated. So, it cannot be used).
 $\mathfrak{C}_k=\mathfrak{C}^+_k\cup \mathfrak{C}_k^-$ and either both  $\mathfrak{C}^+_k$ and $\mathfrak{C}_k^-$ are unbounded, or
$$\mathfrak{C}^+_k\cap \mathfrak{C}_k^- \neq \{( \lambda_k(m, \delta), \theta)\}.$$

We claim that the second case can never occur.

In fact, the uniqueness of IVP guarantees that
$$\mathfrak{C}^+_k\subset \Phi^+_k, \ \ \ \  \mathfrak{C}^-_k\subset \Phi^-_k.
\eqno (3.13)
$$

Suppose on the contrary that $(\eta, z)\in (\mathfrak{C}^+_k\cap \mathfrak{C}_k^-)$ for some  $(\eta, z)\neq (\lambda_k(m, \delta), \theta)$. Then it follows from (3.13) that
 $z=\theta$.
 In this case, $\eta=\lambda_j(m, \delta)$ for some $k\neq j$. Suppose $(\lambda_{m},u_{m})\to (\lambda_j(m, \delta),\theta)$ when $m\to +\infty$ with $(\lambda_{m},u_{m})\in \mathfrak{C}_{k}$. Let $v_{m}=\frac{u_{m}}{\|u_{m}\|_{C^1[\delta, R]}}$, then $v_{m}$ should be a solution of the problem

$$v_{m}=\lambda_{m} \mathcal{L}_\delta(v_{m})+\frac{\mathcal{H}_\delta(\lambda,u_m)}{\|u_{m}\|_{C^1[\delta, R]}}
\eqno(3.14)
$$
This together with the compactness of $\mathcal{L}_\delta$ and $\mathcal{H}_\delta$ imply that for some convenient subsequence $v_{m}\to v_{0}\neq0$ as $m\to +\infty$. Now $v_{0}$ verifies the equation
$$-(r^{N-1}v_{0}')'=\lambda_{j} r^{N-1}m(r)v_{0}$$
and $\|v_{0}\|=1$. Hence $v_{0}\in S_{j, \delta}$. Since $ S_{j, \delta}$ is an open in $E_\delta$, and as a consequence for some $m$ large enough, $v_{m}\in S_{j, \delta}$, and this is a contradiction.

 Therefore, both  $\mathfrak{C}^+_k$ and $\mathfrak{C}_k^-$ are unbounded.

\vskip 3mm

Take
$$\zeta^+:=\mathfrak{C}^+, \ \ \ \ \ \ \zeta^-:=\mathfrak{C}^-.
$$
Obviously, (a) is true.

 (b) can be deduced from the fact that
 $$\sup\{||u'||_{C[\delta, R]}:\, (\lambda, u)\in \zeta^\nu\}\leq 1, \ \ \ \
 \sup\{||u||_{C[\delta, R]}:\, (\lambda, u)\in \zeta^\nu\}\leq R-\delta.
 $$

 (c) Let
 $$\lambda_*^\nu:=\inf\{\lambda: \, (\lambda, u)\in \zeta^\nu\}.
 $$
We claim that $\lambda^\nu_*\in (0, \infty)$.

Suppose on the contrary that  $\lambda^\nu_*=0$. Then there exists a sequence $\{(\mu_n, u_n)\}\subset \zeta^\nu$ satisfying
$u_n>0$, and
$$\lim_{n\to \infty}(\mu_n, u_n)=(0, u^*) \ \ \ \ \ \  \ \text{in}\ \mathbb{R}\times X_\delta
$$
for some $u^*\geq 0$.  Then it follows from
$$-(r^{N-1}\phi_1(u_n'))'=\mu_n r^{N-1} f(r, u_n),\ \ \ \ \ \  u_n'(\delta) = 0 = u_n(R)
$$
that,  after taking a subsequence
and relabeling, if necessary,
$$u_n\to 0.$$
On the other hand,
$$\left\{\begin{array}{ll}
-(r^{N-1}u_n')'=\mu_n r^{N-1}f(r,u_n)h(u_n')-(N-1)r^{N-2}u_n'^{3},\ \ \ \ \ r\in (\delta,R),\\
u_n'(\delta)=0=u_n(R).  \\
\end{array}
\right.
$$
Setting, for
all $n$, $v_n = u_n/||u_n||_{C[\delta,R]}$, we have that
$$\left\{\begin{array}{ll}
-(r^{N-1}v_n')'=\mu_n r^{N-1}\frac{f(r,u_n)}{u_n}h(u_n')v_n-(N-1)r^{N-2}u_n'^{2}\,v_n',\ \ \  \ r\in (\delta,R),\\
v_n'(\delta)=0=v_n(R).  \\
\end{array}
\right.
\eqno(3.15)$$
Notice that
$$r^{N-1}\phi_1(u_n'(r))=-\mu_n\int^r_\delta \tau^{N-1}f(\tau, u_n(\tau))d\tau, \ \ \ \ \ r\in [\delta, R].
$$
This together with $f(r,0)=0$ for $r\in [\delta, R]$ imply that
$$\lim_{n\to\infty}\, ||u'_n||_{C[\delta, R]}=0.
$$
Combining this with (3.15) and the facts $f_0=m(r), \ u_n\to 0$ and $\lim_{n\to\infty} h(u_n')=1$, it concludes that $\mu_n\to \lambda_k(m, \delta)$. This is a contradiction.
  \hfill{$\Box$}

\vskip 5mm

In the following,  we will deal with the case that $f_0=\infty$.

\vskip 2mm

Define $f^{[n]}:[\delta,R]\times \mathbb{R}\to \mathbb{R}$ as follows
 $$f^{[n]}(r,s)=\left\{\begin{array}{lll}
 nf(r,\frac{1}{n})s,\ \ \ \ \  \ &|s|\in [0,\frac{1}{n}],\\
 f(r,s),\ \ \ \ \ \  &  |s|\in (\frac{1}{n},\infty).\\
\end{array}
\right.
$$ Then $f^{[n]}$ is continuous and satisfies (A1) and
$$(f^{[n]})_0=nf(r,\frac{1}{n})=f(r,\frac{1}{n})/(1/n)=:m^{[n]}(r) \ \ \ \text{uniformly for}\ r\in[\delta,R].
$$

Now, let us consider the auxiliary family of the problems
$$\left\{\begin{array}{ll}
-(r^{N-1}u')'=\lambda r^{N-1}f^{[n]}(r,u)h(u')-(N-1)r^{N-2}u'^{3},\ \ \  \ r\in (\delta,R),\\
u'(\delta)=0=u(R).  \\
\end{array}
\right.
\eqno(3.16)
$$
From the definition of $f^{[n]}$, it follows that for $r\in [\delta,R]$ and every $u\in
\mathbb{R}$,
$$ f^{[n]} (r,s)=(m^{[n]}(r)+\xi^{[n]}(r,s))s,$$
where $\xi^{[n]}:[\delta,R]\times \mathbb{R}\to \mathbb{R}$ is continuous and
$$\lim\limits_{s\to0}\xi^{[n]}(r,s)=0\ \ \quad \text{uniformly for}
\ r\in[\delta,R].
$$
Let us set, for convenience, $k(v)=h(v)-1$
for $v\in\mathbb{R}$. We have
$$\lim\limits_{v\to 0}\frac{k(v)}{v}=0.
\eqno(3.17)
$$
Define the operator $\mathcal{H}_\delta^{[n]}:\mathbb{R}\times
E_\delta\to E_\delta$ by
$$\mathcal{H}_\delta^{[n]}(\lambda,u)=\mathcal{G}_\delta\Big(\lambda\big(\xi^{[n]}
(\cdot,u)+[m^{[n]}+\xi^{[n]}(\cdot,u)]k(u')\big)u-\gamma(\cdot)u'^{3}\Big).
$$
Clearly,
$\mathcal{H}_\delta^{[n]}$ is completely continuous  and by (3.16) and
(3.17), it follows that
$$\lim\limits_{\|u\|_{C^1[\delta,R]}\to 0}\frac{\|\mathcal{H}_\delta^{[n]}(\lambda,u)\|_{C^1[\delta,R]}}{\|u\|_{C^1[\delta,R]}}=0,$$
uniformly with respect to $\lambda$ varying in bounded intervals.
Observe that, for any $\lambda$, the couple $(\lambda,
u)\in\mathbb{R}\times E_\delta$ is a solution of the
equation
$$u=\lambda \mathcal{L}_\delta^{[n]}(u)+\mathcal{H}_\delta^{[n]}(\lambda, u)
\eqno(3.18)
$$
if and only if $u$ is a solution of (3.16), here
$\mathcal{L}_\delta^{[n]}:X_\delta\to E_\delta$ be defined by
$\mathcal{L}_\delta^{[n]}(u)=\mathcal{G}_\delta(m^{[n]}u)$.

Let $\Sigma_\delta^{[n]}\subset\mathbb{R}\times E_\delta$ be the closure
of the set of all nontrivial solutions $(\lambda, u)$ of (3.18) with
$\lambda>0$. Note that the set $\{u\in E_\delta | (\lambda,u)\in \Sigma_\delta^{[n]}\}$ is bounded in $E_\delta$.

\vskip3mm

\noindent{\bf Remark 3.2.} Note that from the compactness of the
embedding $E_\delta\hookrightarrow X_\delta$, it concluded that
$\mathfrak{C}_+^{[n]}$ is also  an unbounded connected component in
$[0,\infty)\times X_\delta$.

 \vskip 3mm

\noindent{\bf Proof of Theorem 1.2 with $\delta\in (0, R)$.}\  Similar to the proof of Theorem 1.1 with $\delta\in (0, R)$,  for each fixed $n$ and $\nu\in\{+,-\}$, there exists an
unbounded component $\mathfrak{C}_{\nu,k}^{[n]}\subset\Sigma_\delta^{[n]}$ of solutions of (3.18) joining $(\lambda_k(m^{[n]},\delta),0)\in \mathfrak{C}_{\nu,k}^{[n]}$ to
infinity  in $[0,\infty)\times \mathcal{S}_{k,\delta}$. Moreover,
$(\lambda_k(m^{[n]},\delta),0)\in \mathfrak{C}_{\nu,k}^{[n]}$ is the a
bifurcation point of (3.18) lying on a trivial solution line $u\equiv
0$ and the component $\mathfrak{C}_{\nu,k}^{[n]}\subset\Phi^\nu_k$ joins the infinity in the direction
of $\lambda$ since $u$ is bounded.

It is not difficult to verify
that $\mathfrak{C}^{[n]}_{\nu,k}$ satisfies all conditions in Lemma 2.1 and
consequently $\limsup\limits_{n\to\infty} \mathfrak{C}^{[n]}_{\nu,k}$ contains a
component $\mathfrak{C}_{\nu,k}$  which is unbounded.

From (A3), it follows that for $r\in [\delta,R]$,
$$\lim\limits_{n\to\infty}\frac{f^{[n]} (r,u)}{u}=
\lim\limits_{n\to\infty}\frac{f(r,\frac{1}{n})}{1/ n}=\infty,$$
and consequently,
$$\lim\limits_{n\to\infty}\lambda_k(m^{[n]}, \delta)=0.
\eqno (3.19)
$$
Thus, from (3.19), we have that the component
$\mathfrak{C}_{\nu,k}$ joins $(0,\theta)$ with infinity in the direction of $\lambda$ in $[0,\infty)\times
S_k^\nu$.

We claim that
$$(\mathfrak{C}_{\nu,k}\setminus\{(0, \theta)\})\subset (0, \infty)\times S^\nu_k.
\eqno (3.20)
$$

Suppose on the contrary that there exists a sequence $\{(\mu_n, u_n)\}\subset \mathfrak{C}_{\nu,k}$ satisfying
$$\lim_{n\to \infty}(\mu_n, u_n)=(\mu*, \theta) \ \ \ \ \ \  \ \text{in}\ \mathbb{R}\times X_\delta
$$
for some $\mu*>0$.  Then
$$\left\{\begin{array}{ll}
-(r^{N-1}u_n')'=\mu_n r^{N-1}f^{[n]}(r,u_n)h(u_n')-(N-1)r^{N-2}u_n'^{3},\ \ \ \ \ r\in (\delta,R),\\
u_n'(\delta)=0=u_n(R).  \\
\end{array}
\right.
$$
Setting, for
all $n$, $v_n = u_n/||u_n||_{C[\delta,R]}$, we have that
$$\left\{\begin{array}{ll}
-(r^{N-1}v_n')'=\mu_n r^{N-1}\frac{f^{[n]}(r,u_n)}{u_n}h(u_n')v_n-(N-1)r^{N-2}u_n'^{2}\,v_n',\ \ \  \ r\in (\delta,R),\\
v_n'(\delta)=0=v_n(R).  \\
\end{array}
\right.
\eqno (3.21)
$$
Notice that
$$r^{N-1}\phi_1(u_n'(r))=-\mu_n\int^r_0 \tau^{N-1}f^{[n]}(\tau, u_n(\tau))d\tau, \ \ \ \ \ r\in [0, R].
\eqno (3.22)
$$
This together with $f^{[n]}(r,0)=0$ for $r\in [\delta, R]$ imply that
$$\lim_{n\to\infty}\, ||u'_n||_{C[0, R]}=0.
\eqno (3.23)
$$
Combining this with (3.21) and the facts $f_0=\infty$ and $\lim_{n\to\infty} h(u_n')=1$, it concludes that
 $\mu^*=0$. This is a contradiction.

 Therefore, (3.20) holds.
\hfill$\Box$

 \vskip 3mm

\section{Radial solutions for the prescribed mean curvature problem in a ball}

In this section, we shall deal with $(1.6)_\delta$ with $\delta=0$.

\vskip 3mm

Let
$$g_n(r, s)
=\left\{
\aligned
0,\qquad \ \ \ \ \ \ \ \ \ \ \ \ & (r, s)\in (0, \frac 1n]\times (-\alpha, \alpha),\\
f(r-\frac 1n,s), \ \ \ \ \ \  &(r, s)\in (\frac 1n, R)\times (-\alpha, \alpha).\\
\endaligned
\right.
\eqno (4.1)
$$
In the following, we shall use the solutions of the family of problems
$$
\left\{\aligned
&-(r^{N-1}\phi_{1}(u'))'=\lambda r^{N-1}g_n(r,u),\ \ \ \ \ r\in (\frac 1n, R),\\
&\ u'(1/n)=0=u(R),\\
\endaligned
\right.
\eqno (4.2)_n
$$
to construct the radial solutions of the prescribed mean curvature problem in a ball
$$\mathcal{M}v+\lambda f(|x|, v)=0 \ \ ~~~\text{in} ~~\mathcal{B}(R), ~~~ v=0 ~~~\text{on} ~~\partial\mathcal{B}(R).
\eqno (4.3)
$$

To find a radial solution of (4.3), it is enough to find a solution of the problem
$$
\left\{\aligned
&-(r^{N-1}\phi_{1}(u'))'=\lambda r^{N-1}f(r,u),\\
&\ u'(0)=0=u(R).\\
\endaligned
\right.
\eqno (4.4)
$$

\vskip 3mm

\ For given $n\in \mathbb{N}$ , let
$(\lambda,u)$ be a solution of $(4.2)_n$.
 For each $n$,
define a function $y_n:[0, R]\to [0, \infty)$ by
$$y_n(r)
=\left\{\aligned
u(r),  \ \ \ \ \ \ \ \ \ \ \ &\frac 1n\leq r\leq R,\\
u(\frac 1n),  \ \ \ \ \ \ \ \ \ \ &0\leq r\leq \frac 1n.
\\
\endaligned
\right.
\eqno (4.5)_n$$
Then
$$y_n\in \{w\in C^2[0, R]:\,w'(0)=w(R)=0\}.$$
Moreover, $y_n$ is a solution of the problem
$$
\left\{\aligned
&-(r^{N-1}\phi_{1}(u'))'=\lambda r^{N-1}g_n(r,u),\ \ \ \ \ r\in (0, R)\\
&\ u'(0)=0=u(R),\\
\endaligned
\right.
\eqno (4.6)_n
$$
i.e. $y_n$ is a solution of the problem
$$
\left\{\aligned
&-(r^{N-1}u'(r))'+ (N-1)r^{N-2} [u'(r)]^3=\lambda r^{N-1} g_n(r, u(r)) h(u'(r)),\ \ \ \ r\in (0, R),
\\
& u'(0)=u(R)=0.
\endaligned
\right.
\eqno (4.7)_n
$$
On the other hand, if $(\lambda, y)$ is a solution of $(4.7)_n$, then $(\lambda, y|_{[\frac 1n, R]})$ is a solution of $(4.2)_n$.

\noindent{\bf Lemma 4.1.} \ Let (A1) and (A2) hold. Let
$\hat\lambda:\, \hat\lambda\neq \lambda_k(m,0)$ be given. Then there exists $\hat b>0$, such that
$$||u||_{C[0, R]}\geq \hat b$$
 for any
 solution $(\hat\lambda, u)\in \Phi_{k,0}^\nu$ of $(4.7)_n$. Here $b$ is independent of $n$ and $u$.

\noindent{\bf Proof}. Suppose on the contrary that $(4.7)_{n}, \
n\in \mathbb{N}$, has a sequence of solution $(\hat\lambda,
y_j)\in \Phi_{k,0}^\nu$ with
$$\lim_{j\to\infty}\, ||y_j||_{C[0, R]}=0.
\eqno (4.8)$$
Then
$$
\left\{\aligned
&(r^{N-1}\phi_1(y_j'(r)))'+\hat\lambda r^{N-1} g_n(r, y_j(r))=0,\ \ \ \ r\in (0, R),\\
&y_j'(0)=y_j(R)=0,\\
\endaligned
\right.
\eqno (4.9)
$$
and consequently,
$$r^{N-1}\phi_1(y_j'(r))=-\hat\lambda\int^r_0 \tau^{N-1}g_n(\tau, y_j(\tau))d\tau, \ \ \ \ r\in [0, R].$$
This together with (4.8) and the fact that $g_n(r,0)=0$ for $r\in [0, R]$ imply that
$$\lim_{j\to\infty}\, ||y'_j||_{C[0, R]}=0.
\eqno (4.10)
$$
Recall that (4.9) can be rewritten as
$$
\left\{\aligned
&-(r^{N-1}y_j'(r))'+(N-1)r^{N-2} [y_j'(r)]^3=\hat\lambda r^{N-1} g_n(r, y_j(r)) h(y_j'(r)),\\
& y_j'(0)=y_j(R)=0.
\endaligned
\right.
\eqno (4.11)_n
$$
Setting, for all $j$, $v_j = y_j/||y_j||_{C[0,R]}$ , we have that
$$
\left\{\aligned
&-(r^{N-1}v_j'(r))'+(N-1)r^{N-2} [y_j'(r)]^2v_j'(r)
=\hat\lambda r^{N-1} \frac{g_n(r, y_j(r))}{y_j(r)}v_j(r) h(y_j'(r)),\\
& v_j'(0)=v_j(R)=0.
\endaligned
\right.
\eqno (4.12)_n
$$
Letting $j\to\infty$, it follows from (4.8), (4.10) and $(4.12)_n$ that there exists $w\in C^2[0,R]$ with $||w||_{C[0, R]}=1$, such that
$$
\left\{\aligned
&-(r^{N-1}w'(r))'=\hat \lambda r^{N-1} m(r)w(r),\\
& w'(0)=w(R)=0,
\endaligned
\right.
\eqno (4.13)
$$
which implies that $\hat\lambda=\lambda_k(m, 0)$. However, this contradicts the assumption $\hat\lambda\neq \lambda_k(m, 0)$.
\hfill{$\Box$}

\vskip 3mm
 Using the same argument with obvious changes, we may prove the following
\vskip 3mm
\noindent{\bf Lemma 4.2.} \ Let (A1) and (A3) hold. Let
$\hat\lambda\in (0, \infty)$ be given. Then there exists $\hat b>0$, such that
$$||u||_{C[0, R]}\geq \hat b$$
 for any
solution $(\hat\lambda, u)\in \Phi^\nu_{k,0}$ of $(4.7)_n$.
\hfill{$\Box$}
\vskip 3mm

Now, we are in the position to prove Theorem 1.1-1.2 with $\delta=0$.

\vskip 3mm

\noindent{\bf Proof of Theorem 1.1 with $\delta=0$.} \ For given $n$, let
$\xi_n$ be the component obtained by Theorem 1.1 with $\delta\in (0, R)$ for $(4.2)_n$.
Let
$$\zeta_n:=\{(\lambda, y_n):\, \ y_n\ \text{is determined by} \ u \ \text{via}\ (4.5)_n  \ \text{for}\ (\lambda,u)\in \xi_n \}.
$$
Then $\zeta_n$ is a component in $[0, \infty)\times C^1[0,R]$ which joins $(\lambda_k(m^{[n]}, \frac 1n), \theta)$ with infinity in the direction of $\lambda$ and
$$\sup\{||y||_{C^1[0, R]}:(\lambda, y)\in \xi_n\}<M
\eqno (4.14)
$$
for some constant $M>0$, independent of $y$ and $n$. Here
$$m^{[n]}(r):=m(r-\frac 1n),  \ \ \ \ \ \ \ \ \ \ \ \frac 1n\leq r\leq R,$$
and $\lambda_k(m^{[n]}, \frac 1n)$ is the principal eigenvalue of the linear problem
$$
\left\{\aligned
&-(r^{N-1}u'(r))'=\lambda r^{N-1} m^{[n]}(r)u(r),\ \ \ \ r\in (\frac 1n, R),\\
&u'(\frac 1n)=u(R)=0.
\endaligned
\right.
\eqno (4.15)
$$
Since
$\lim_{n\to\infty}\, \lambda_k(m^{[n]}, \frac 1n)=\lambda_k(m, 0)$, it follows from Lemma 2.1 that there exists a component $\zeta$ in
     \ $\underset{n\to\infty}{\limsup}\; \zeta_n$ \ which joins  $(\lambda_k(m,0), \theta)$ with infinity in the direction of $\lambda$ and
$$\sup\{||y||_{C^1[0, R]}:(\lambda, y)\in \zeta\}\leq M.
\eqno (4.16)
$$
Now, Lemma 4.1 ensures that
$$\zeta\cap\big([0, \infty)\times \{\theta\}\big)=\{(\lambda_k(m, 0),\theta)\}.
$$
 \hfill{$\Box$}

 \noindent{\bf Proof of Theorem 1.2 with $\delta=0$.} \ It is an immediate consequence of Theorem 1.2 with $\delta>0$ and Lemma 4.2.
  \hfill{$\Box$}

\vskip 3mm

The proof of Corollary 1.1 is a direct consequence of Theorem 1.1.

\vskip 3mm

\noindent{\bf Proof of Corollary 1.2.} \ From Theorem 1.2, for each $\nu\in\{+,-\}$ and $n\in \mathbb{N}$, $(1.6)_\delta$ has a solution $u_n^\nu\in S^\nu_{n,\delta}$. We only need to show that
$$\lim_{n\to\infty}\; |u^\nu_n(r)|+|(u_n^\nu)'(r)|=0 \ \ \text{uniformly in}\ r\in [0,R].
$$

\vskip 2mm

We firstly deal with the case that $\delta>0$. In the following, we shall replace $u_n^\nu$ with $u_n$
for fixed $\nu\in \{+,-\}$.

\vskip 2mm

{\it Step 1} We show that
$$\lim_{n\to\infty}||u_n||_{C[\delta, R]}=0.
\eqno (4.17)$$

\vskip 2mm

Let $\tau_1(n), \cdots, \tau_n(n)$ be zeros of $u_n$ in $(0, R]$:
$$\big(\delta=\tau_0(n)<\big)\, \tau_1(n)<\cdots <\tau_{n-1}(n)<\tau_n(n)\; \big(=R\big).
\eqno (4.18)$$
Since $||u_n'||_{C[\delta, R]}<1$, it follows that
     $$||u_n||_{C[\delta, R]}\leq \sup\{|\tau_{j+1}(n)-\tau_j(n)|:\;j=0, \cdots, n-1\}.
  \eqno (4.19)
     $$

\vskip 2mm

Suppose on the contrary that (4.17) is not true. Then there exist a positive constant $\sigma_0$ and a subsequence of solution of $(1.6)_\delta$, $\{(1, u_{n_j})\}\subseteq \{(1, u_{n})\}$, such that
$$||u_{n_j}||_{C[\delta, R]}\geq \sigma_0.
\eqno (4.20)
$$
It follows from (4.20)
 and (4.19) that for $n_j\geq 6$, there exist $t_1(n_j), t_2(n_j) \in \{\tau_0(n_j), \tau_1(n_j), \cdots, $ $\tau_{n_j}(n_j)\}$ satisfying
 $$t_2(n_j)-t_1(n_j)\geq \sigma_0,
 \eqno (4.21)
 $$
 $$||u_{n_j}||_{C[\delta, R]}=\sup\{|u_{n_j}(r)|:\, r\in [t_1(n_j),t_2(n_j)]\},
 \eqno (4.22)
 $$
 $$u_{n_j}(r)\neq 0, \ \ \ \  r\in (t_1(n_j),t_2(n_j)).
 \eqno (4.23)
 $$
 Without loss of generality, we may assume that there exists a closed subinterval $I_1$ with $
 \text{meas}\, I_1\geq \frac {\sigma_0}2$,
  such that (after taking a subsequence and relabeling, if necessary,)
 $$u_{n_j}(r)\neq 0, \ \ \ \ \  r\in I_1, \  j\geq j_0 \ \text{for some} \ j_0\in \mathbb{N}.
 \eqno (4.23)
 $$

Now, applying the facts
 $$u_{n_j}(r)=\int_r^R\phi_1^{-1}\big(\frac 1{s^{N-1}}\int^s_\delta t^{N-1}f(t, u_{n_j}(t)dt)\big)ds
 \eqno (4.24)
 $$
 and $$||u'_{n_j}||_{C[\delta, R]}<1, \ \ \ \ ||u_{n_j}||_{C[\delta, R]}<R,
 $$
 and the standard argument, we deduce that after taking a subsequence
and relabeling, if necessary,
$$u_{n_j}\to u_\diamond\ \ \ \  \text{in} \ C^1[\delta, R],
$$
for some $u_\diamond\in C^1[\delta, R]$. Obviously,
$$u_\diamond(r)\geq \frac {\sigma_0}2, \ \ \ \ \ r\in I_1,
\eqno (4.25)
$$
 and it is easy to check that $u_\diamond$ satisfies
 $$
\left\{\aligned
&-(r^{N-1}u_\diamond'(r))'= r^{N-1}f(r,u_\diamond),\ \ \ \ r\in (\delta, R),\\
&u_\diamond'(\delta)=u_\diamond(R)=0.
\endaligned
\right.
$$

 On the other hand, we may take $s_1(n_j), s_2(n_j) \in \{\tau_1(n_j), \cdots, \tau_{n_j}(n_j)\}$ with  $$s_2(n_j)>s_1(n_j),
 \eqno (4.26)$$
 $$s_2(n_j),s_1(n_j)\not\in [t_1(n_j),t_2(n_j)],
 \eqno (4.27)$$
 $$u_{n_j}(r)\neq 0, \ \ \ \  r\in (s_1(n_j),s_2(n_j)),
 \eqno (4.28)
 $$
  $$s_2(n_j)-s_1(n_j)\to 0, \ \ \ \ \ j\to\infty.
 \eqno (4.29)
 $$
 For each $n_j\geq 6$, $u_{n_j}(s_1(n_j))=u_{n_j}(s_2(n_j))=0$ yields
 $$u'_{n_j}(x(n_j))=0, \ \ \ \ \text{for some}\ x(n_j)\in((s_1(n_j),\; s_2(n_j)).
 \eqno(4.30)
 $$
After taking a subsequence
and relabeling, if necessary, we may assume
$$x(n_j)\to x_0, \ \ \ \ \ j\to \infty.$$
Combining this with (4.30), (4.29) and using (4.24), it follows that
$$u_\diamond(x_0)=u'_\diamond(x_0)=0\ \ \ \ \text{for some}\ x_0\in [\delta,R]
\eqno (4.31)
$$
which means that $u_\diamond\equiv 0$ in $[\delta, R]$ (see Lemma 2.5 where $x_0>0$ is needed). However, this contradicts (4.25).

Therefore, (4.17) is true.

 \vskip 2mm

 {\it Step 2} We show that
$$\lim_{n\to\infty}||u'_n||_{C[\delta, R]}=0.
\eqno (4.32)
$$

  In fact, it is an immediate consequence of (4.17) and the relation
  $$u_n'(r)=-\phi_1^{-1}\big(\frac 1{r^{N-1}}\int^r_0 t^{N-1}f(t, u_n(t)dt)\big).
  $$

  \vskip 4mm

  Next, we consider the case $\delta=0$.

  \vskip 2mm
  Let
  $$g_m(r, s)
=\left\{
\aligned
0,\qquad \ \ \ \ \ \ \ \ \ \ \ \ & (r, s)\in (0, \frac 1m]\times (-\alpha, \alpha),\\
f(r-\frac 1m,s), \ \ \ \ \ \  &(r, s)\in (\frac 1m, R)\times (-\alpha, \alpha).\\
\endaligned
\right.
$$
Let $u^\nu_m\in S^\nu_{n, \frac 1m}$. Then
  $$
\left\{\aligned
&-(r^{N-1}\phi_{1}((u^\nu_m)'))'= r^{N-1}g_m(r,u^\nu_m),\ \ \ \ \ r\in (\frac 1m, R),\\
&\ (u^\nu_m)'(1/m)=0=u^\nu_m(R).\\
\endaligned
\right.
$$
For each $n$ and $m$,
define a function $y^\nu_m\in S^\nu_{n, 0}$ by
$$y^\nu_m(r)
=\left\{\aligned
u_m(r),  \ \ \ \ \ \ \ \ \ \ \ &\frac 1m\leq r\leq R,\\
u_m(\frac 1m),  \ \ \ \ \ \ \ \ \ \ &0\leq r\leq \frac 1m.
\\
\endaligned
\right.
$$
 By the same argument used in proof of Theorem 1.1, with obvious changes, we may use $y^\nu_m$ to construct two solutions, $z_n^\nu\in S^\nu_{n, 0}, \  \nu\in\{+, -\}$, for $(1.5)_{\lambda=1}$.
 Obviously, (4.17) and (4.32) imply that
  $$\lim_{n\to\infty}(||z^\nu_n||_{C[0, R]}+||(z^\nu)'_n||_{C[0, R]})=0.
$$

\hfill{$\Box$}



\newpage

\centerline {\bf REFERENCES}\vskip5mm\baselineskip 20pt
\begin{description}

\item{[1]}\ L. J. A\'{i}as, B. Palmer, On the Gaussian curvature of maximal surfaces and the Calabi-Bernstein theorem, Bull. London Math. Soc. 33 (2001) 454-458.

\item{[2]} A. Ambrosetti, J. Garcia-Azorero and I. Pekal, Quasilinear equations with a multiple bifurcation,
Differential Integral Equations 10 (1997) 37-50.

\item{[3]}\ V. Anuradha, D. D. Hai,  R. Shivaji,
Existence results for superlinear semipositone BVP's,
Proc. Amer. Math. Soc. 124(3) (1996)  757-763.

\item{[4]}\ R. Bartnik, L. Simon, Spacelike hypersurfaces with prescribed boundary values and mean curvature, Comm. Math.
Phys. 87 (1982-1983) 131-152.

\item{[5]}\  C. Bereanu, P. Jebelean, J. Mawhin, Radial solutions for Neumann problems involving mean curvature operators in
Euclidean and Minkowski spaces, Math. Nachr. 283 (2010) 379-391.

\item{[6]}\ C. Bereanu, P. Jebelean, J. Mawhin, Radial solutions for some nonlinear problems involving mean curvature operators
in Euclidean and Minkowski spaces, Proc. Amer. Math. Soc. 137 (2009) 171-178.

\item{[7]}\ C. Bereanu, P. Jebelean, P. J. Torres, Positive radial
 solutions for  Dirichlet problems with mean curvature operators in
 Minkowski space, J. Funct. Anal. 264 (2013) 270-287.

\item{[8]}\ C. Bereanu, P. Jebelean, P. J. Torres, Multiple positive
 radial solutions for a Dirichlet problem involving the mean curvature
 operator in Minkowski space, J. Funct. Anal. 265(4) (2013) 644-659.

\item{[9]}\ M. F. Bidaut-V\'{e}ron, A. Ratto, Spacelike graphs with prescribed mean curvature, Differential Integral Equations 10
(1997) 1003-1017.

\item{[10]}\ A. Capietto, W. Dambrosio, F. Zanolin, Infinitely many radial solutions to a boundary value problem in a ball, Ann. Mat. Pura Appl. (4) 179 (2001) 159-188.

\item{[11]}\  A. Castro, A. Kurepa, Infinitely many radially symmetric solutions to a superlinear Dirichlet problem in a ball, Proc. Amer. Math. Soc. 101 (1987)  57-64.

\item{[12]}\ S.-Y. Cheng, S.-T. Yau, Maximal spacelike hypersurfaces in the Lorentz-Minkowski spaces, Ann. of Math. 104 (1976) 407-419.

\item{[13]}    Y. Cheng, On the existence of radial solutions of a nonlinear elliptic equation on the unit ball, Nonlinear
Anal. 24 (1995) 287-307.

\item{[14]}\ I. Coelho, C. Corsato, F. Obersnel, P. Omari,  Positive solutions of the Dirichlet problem for the one-dimensional Minkowski-curvature equation, Adv. Nonlinear Stud. 12(3) (2012)  621-638.

\item{[15]}\ G. Dai, \ R. Ma, \ Unilateral global bifurcation phenomena and nodal solutions for $p$-Laplacian, J. Differential Equations 252(3) (2012) 2448-2468.

\item{[16]}\ E. N. Dancer, Bifurcation from simple eigenvalues and eigenvalues of geometric multiplicity one, Bull. Lond. Math. Soc. 34 (2002) 533-538.

\item{[17]}\  E. N. Dancer On the structure of solutions of non-linear eigenvalue problems, Indiana Univ. Math. J. 23 (1973/74) 1067-1076.

\item{[18]}\ M. J. Esteban, Multiple solutions of semilinear elliptic problems in a ball, J. Differential Equations
     57 (1985) 112-137.

\item{[19]} M. Garc\'{i}a-Huidobro, R. Man\'{a}sevich,  F. Zanolin, Strongly nonlinear second-order ODEs with rapidly growing terms, J. Math. Anal. Appl. 202(1) (1996)  1-26.

\item{[20]} M. Grillakis, Existence of nodal solutions of semilinear equations in $\mathbb{R}^N$, J. Differential Equations 85 (1990) 367-400.

\item{[21]} Z. Guo, Boundary value problems for a class of quasilinear ordinary differential equations, Differential
     Integral Equations  6 (1993) 705-719.

\item{[22]}\ M. A. Krasnosel'skii,  A. I. Perov,  A. I. Povolotskiy, P. P. Zabreiko,
           Plane vector fields, Translated by Scripta Technica, Ltd. Academic Press, New York-London 1966.

\item{[23]}\ R. L\'{o}pez, Stationary surfaces in Lorentz-Minkowski space, Proc. Roy. Soc. Edinburgh Sect. A 138 (2008) 1067-1096.

\item{[24]}\ J. L\'{o}pez-G\'{o}mez, Spectral Theory and Nonlinear Functional Analysis, in: Research Notes in Mathematics, vol. 426, Chapman \& Hall/CRC, Boca Raton,
    Florida, 2001.

\item{[25]} J. L\'{o}pez-G\'{o}mez,  C. Mora-Corral, Algebraic Multiplicity of Eigenvalues of Linear Operators, Oper. Theory Adv. Appl. vol. 177,
    Birkh\"{a}user/Springer, Basel, Boston/Berlin, 2007.

\item{[26]}\ R. Ma, Y. An, Global structure of positive solutions for superlinear second order $m$-point boundary value problems, Topol. Methods Nonlinear Anal. 34(2) (2009) 279-290.

\item{[27]}\ R. Ma, Y. An,  Global structure of positive solutions for nonlocal boundary value problems involving integral conditions, Nonlinear Anal. 71(10) (2009)  4364-4376.

\item{[28]}\   R. Ma, C. Gao, Bifurcation of positive solutions of a nonlinear discrete
      fourth-order boundary value problem, Z. Angew. Math. Phys. 64 (2013) 493-506.

\item{[29]}\ J. Mawhin, Radial solution of Neumann problem for periodic perturbations of the mean extrinsic curvature operator,
     Milan J. Math. 79 (2011) 95-112.

\item{[30]} \ F. I. Njoku, P. Omari, F. Zanolin,
       Multiplicity of positive radial solutions of a quasilinear elliptic problem in a ball,
       Adv. Differential Equations, 5 (2000) 1545-1570.

\item{[31]}\ P. Omari, F. Zanolin,
     Infinitely many solutions of a quasilinear elliptic problem with an oscillatory
      potential, Comm. Partial Differential Equations 21 (1996) 721-733.

\item{[32]}\ I. Peral, Multiplicity of solutions for the $p$-Laplacian, ICTP SMR 990/1, 1997.

\item{[33]}\ P. H. Rabinowitz, Some global results for nonlinear eigenvalue problems, J. Funct. Anal. 7 (1971) 487-513.

\item{[34]}\  B. L. Shekhter,
     On existence and zeros of solutions of a nonlinear two-point boundary value problem,
     J. Math. Anal. Appl. 97 (1983) 1-20.

\item{[35]}\ A. E. Treibergs, Entire spacelike hypersurfaces of constant mean curvature in Minkowski space, Invent. Math. 66 (1982) 39-56.

\item{[36]}\ G. T. Whyburn, Topological Analysis, Princeton University Press, Princeton, 1958.

\item{[37]} M. Willem, Minimax theorems, Birkh\"{a}user, Boston, 1996.

\end{description}
\end{document}